\documentclass[review]{elsarticle}
\newtheorem{thm}{Theorem}

\newtheorem{lem}{Lemma}

\newtheorem{prop}[thm]{Proposition}
\newtheorem{Definition}{Definition}
\newtheorem{rem}{Remark}

\newcommand{\B}{\mathcal{B}}

\newcommand{\A}{\mathcal{A}}
\newcommand{\I}{\mathcal{I}}
\newcommand{\D}{\mathcal{D}}
\newcommand{\F}{\mathcal{F}}

\newcommand{\C}{\mathcal{C}}
\newcommand{\M}{\mathcal{M}}
\newcommand{\N}{\mathcal{N}}
\newcommand{\Z}{\mathcal{Z}}
\newcommand{\E}{\mathcal{E}}

\newcommand{\norm}[1]{\left\Vert#1\right\Vert}

\usepackage{graphicx}
\usepackage{graphics}
\usepackage{color}
\usepackage{dsfont}
\usepackage{mathrsfs}
\usepackage{amsfonts}
\usepackage{multirow}
\usepackage{graphicx}   
\usepackage{epsfig}
\usepackage{multicol}
\usepackage{amssymb}
\usepackage{amsmath}
\usepackage{anysize}
\usepackage{graphicx}
\usepackage{caption}
\usepackage{subfig}
\usepackage{tcolorbox}
\begin{document}
\begin{frontmatter}
\title {Fast and flexible preconditioners for solving multilinear systems}
\author {Eisa Khosravi Dehdezi}
\ead{esakhosravidehdezi@gmail.com}

\author {Saeed Karimi \corref{mycorrespondingauthor}}
\address{Department of Mathematics, Persian Gulf University, Bushehr, Iran}
\cortext[mycorrespondingauthor]{Corresponding author}
\ead{karimi@pgu.ac.ir}

\begin{abstract}
This paper investigates a type of fast and flexible preconditioners to solve multilinear system $\A\textbf{x}^{m-1}=\textbf{b}$ with $\M$-tensor $\A$ and obtains some important convergent theorems about preconditioned Jacobi, Gauss-Seidel and SOR type iterative methods. The main results theoretically prove that the preconditioners can accelerate the convergence of iterations. Numerical examples are presented to reverify the efficiency of the proposed preconditioned methods.
\end{abstract}
\end{frontmatter}
\section{introduction}
In recent decades, Tensors or hypermatrices have been applied in many types of research and application areas such as data analysis, psychometrics, chemometrics, image processing, graph theory, Markov chains, hypergraphs, etc. \cite{R30}. Tensor equations (or multilinear systems \cite{R5}) involving the Einstein product have been discussed in \cite{A3}, which has many applications in continuum physics, engineering, isotropic and anisotropic elastic models \cite{R21}. Wang and Xu presented some iterative methods for solving several kinds of tensor equations in \cite{R38}, Huang and Ma, in \cite{R15}, proposed the Krylov subspace methods to solve a class of tensor equations. In \cite{R19}, Khosravi Dehdezi and Karimi proposed the extended conjugate gradient squared and conjugate residual squared methods for solving the generalized coupled Sylvester tensor equations
\[ \sum_{j=1}^{n}
 \mathcal{X}_j\times_1A_{ij1}\times_2A_{ij2}\times... \times_dA_{ijd}=\mathcal{C}_{i},\hspace{0.5cm}i=1,2,...,n,\]
where the matrices $A_{ijl}\in\Bbb{C}^{n_{ijl}\times n_{ijl}}$ ($i,j=1,2,...,n$ and $l=1,2,...,d$), tensors $\mathcal{C}_{i}\in\Bbb{C}^{n_{i1}\times ...\times n_{id}}(i=1,2,...,n)$ are known, tensors $\mathcal{X}_j\in\Bbb{C}^{n_{j1}\times...\times n_{jd}}(j=1,2,...,n)$ are unknown and $\times_j(j=1,2,..,n)$ is the $j$-mode product. Also they proposed a fast and efficient Newton-Shultz-type iterative method for computing inverse and Moore-Penrose inverse of tensors in \cite{R78}.\\
\indent
Very recently years, solving the following multilinear system has become a hot topic because of several applications such as data analysis, engineering and scientific computing \cite{A1,A2,A3}:
\begin{equation}\label{eq1}
\A\textbf{x}^{m-1}=\textbf{b}
\end{equation}
where $\A=(a_{i_1i_2…i_m})$ is an $m$ order $n$-dimensional tensor, $\textbf{x}$ and $\textbf{b}$ are vectors in $\Bbb{C}^n$. The $n$ dimensional vector $\A\textbf{x}^{m-1}$ is defined as \cite{A4}:
\begin{equation}\label{eq2}
(\A\textbf{x}^{m-1})_i=\sum_{i_2=1}^{n}...\sum_{i_m=1}^{n}
a_{ii_2...i_m}x_{i_2}...x_{i_m},~~~i=1,2,...,n,
\end{equation}
and $x_i$ denotes the $i$-th component of $\textbf{x}$.\\
\indent
Many theoretical analyses and algorithms for solving \eqref{eq1} were also studied. Qi in \cite{A4} considered an $m$ order $n$-dimensional supersymmetric tensor and showed that when $m$ is even it has exactly $n(m-1)^{n-1}$ eigenvalues, and the number of its E-eigenvalues is strictly less than $n(m-1)^{n-1}$ when $m\geq4$. Ding and Wei in \cite{A3} proved that a nonsingular $\M$-equation with a positive right-hand side always has a unique positive solution. Also, they applied the $\M$-equations to some nonlinear differential equations and the inverse iteration for spectral radii of nonnegative tensors. In \cite{A6}, Han proposed a homotopy method for ﬁnding the unique positive solution to a multilinear system with a nonsingular $\M$-tensor and a positive right side vector. Li et al., in \cite{A10} extended the Jacobi, Gauss‐Seidel and successive over‐relaxation (SOR) iterative methods to solve the tensor equation $\A\textbf{x}^{m-1}=\textbf{b}$, where $\A$ is an $m$ order $n$-dimensional symmetric tensor. Under appropriate conditions, they showed that the proposed methods were globally convergent and locally $r$‐linearly convergent. In \cite{A7}, He et al. proved that solving multilinear systems with $\M$-tensors is equivalent to solving nonlinear systems of equations where the involving functions are P-functions. Based on this result, they proposed a Newton-type method to solve multilinear systems with $\M$-tensors. For a multilinear system with a nonsingular $\M$-tensor and a positive right side vector, they showed that the sequence generated by the method converges to the unique solution of the multilinear system and the convergence rate is quadratic. For solving the multilinear systems, Liang et al. in \cite{A11}, transformed equivalently the tensor equation into a consensus constrained optimization problem, and then proposed an ADMM type method for it. Also, they showed that each limit point of the sequences generated by the method satisfied the Karush-Kuhn-Tucker conditions.
Liu et al., in \cite{A8}, introduced the variant tensor splittings, and presented some equivalent conditions for a strong $\M$-tensor based on the tensor splitting. Also, the existence and unique conditions of the solution for multi-linear systems were given.
Besides, they proposed some tensor splitting algorithms for solving multi-linear systems with coefficient tensor being a strong $\M$-tensor. As an application, a tensor splitting algorithm for solving the multi-linear model of higher-order Markov chains was proposed. Li et al., in \cite{A9} firstly derived a necessary and sufficient condition for an $\M$-tensor equation to have nonnegative solutions. Secondly, developed a monotone iterative method to find a nonnegative solution to an $\M$-tensor
equation. Under appropriate conditions, they showed that the sequence of iterates generated by the method converges to a nonnegative solution of the $\M$-tensor equation monotonically and linearly. Bai et al. in \cite{A4.5}
proposed an algorithm that always preserves the nonnegativity of solutions of the multilinear system under consideration involves a nonsingular $\M$-tensor and a nonnegative right-hand side vector. Also, they proved that the sequence generated by the proposed algorithm is a nonnegative componentwise nonincreasing sequence and converges to a nonnegative
solution of the multilinear system. Cui et al. in \cite{A5} intended to solve the multi-linear system by the preconditioned iterative method based on tensor splitting. For this purpose, they proposed the preconditioner $I+S_{\max}$. Lv and Ma in \cite{A12} proposed a Levenberg-Marquardt (LM) method for solving tensor equations with semi-symmetric coefficient tensor and proved its global convergence and local quadratic convergence under the local error bound condition, which is weaker than non-singularity. As an application, they solved the H-eigenvalue of real semi-symmetric tensor by the LM method. Wang et al., in \cite{A13} proposed continuous-time neural network and modified continuous-time neural networks for solving a multi-linear system with $\M$-tensors. They proved that the presented neural networks are stable in the sense of Lyapunov stability theory. For solving the multilinear system $\A\textbf{x}^{m-1}=\textbf{b}$, where $\A$ is a symmetric $\M$-tensor, Xie et al. in \cite{A14} proposed some tensor methods based on the rank-1 approximation of the coefﬁcient tensor. Li et al. in \cite{A15}, considered tensor equations of 3 order whose solutions
are the intersection of a group of quadrics from a geometric point of view. Inspired by the method of alternating projections for set intersection problems, they developed a hybrid alternating projection algorithm for solving these tensor equations. The local linear convergence of the alternating projection method was established under suitable conditions. Liu et al. in \cite{A16}, presented a preconditioned SOR method for solving the multilinear systems whose coefficient tensor is an $\M$-tensor. Also, the corresponding comparison for spectral radii of iterative tensors was given. It is known that the preconditioning technique plays an important role in solving multilinear systems. In particular, when the coefficient tensor is an $\M$-tensor, there is little research on these techniques so far. By this motivation, we establish some effective preconditioners and give a theoretical analysis.\\
\indent
The rest of this paper is organized as follows. Section 2 is preliminary in which we introduce some related definitions and lemmas. In Section 3, new fast and flexible type preconditioners are proposed, and the corresponding theoretical analysis is given. In Section 4, numerical examples are given to show the efficiency of the proposed preconditioned iterative methods. Section 5 is the concluding remark and the final section is the future researches.
\section{Preliminaries}
In this section, we introduce some definitions, notations, and related properties which will be used in the following.\\
\indent
Let \textbf{0}, $O$ and $\mathcal O$ denote for null vector, null matrix and null tensor, respectively. Let $\A$ and $\B$ be a tensor (vector or matrix) with the same sign. The order $\A\geq\B(>\B)$ means that each element of $\A$ is no less than (larger than) corresponding one of $\B$.\\
\indent
A tensor $\A$ consists of $n_1\times...\times n_m$ elements in the complex field $\Bbb C$:
\[\A=(a_{i_1i_2…i_m}),~~a_{i_1i_2…i_m}\in\Bbb C,~~1\leq i_j\leq n_j,j=1,...,m.\]
When $m=2$, $\A$ is an $n_1\times n_2$ matrix. If $n_1=...=n_m=n$, $\A$ is called an $m$ order $n$-dimensional tensor. By $\Bbb C^{n_1\times...\times n_m}$ we denote all $m$ order tensors consist of $n_1\times...\times n_m$ entries and by $\Bbb C^{[m,n]}$ we denote the set of all $m$ order $n$-dimensional tensors. When $m=1$, $\Bbb C^{[1,n]}$ is simplified as $\Bbb C^n$ , which is the set of all $n$-dimension complex vectors. Similarly, the above notions can be used to the real number field $\Bbb R$.\\
\indent
Let $\A\in\Bbb R^{[m,n]}$. If each entry of $\A$ is nonnegative, then $\A$ is called a nonnegative tensor. The set of all $m$ order $n$-dimensional nonnegative tensors is denoted by $\Bbb R_+^{[m,n]}$. The $m$ order $n$-dimensional identity tensor, denoted by $\mathcal I_m=(\delta_{i_1i_2...i_m})\in\Bbb R^{[m,n]}$, is the tensor with entries:
\[\delta_{i_1i_2...i_m}=\left\lbrace\begin{array}{l}
                               1, ~~i_1=i_2=...=i_m \\
                               0, ~~otherwise.
                               \end{array}\right.\]
When $m=2$, the identity tensor reduces to identity matrix of size $n\times n$, denoted by $I$.
\begin{Definition}\label{d01}\cite{A20}
$\A\in\Bbb C^{[m,n]}$ is called a reducible tensor if there exists a nonempty proper index subset $\mathbb I\subseteq\{1,2,...,n\}$ such that
\[a_{i_1i_2...i_m}=0,~ \forall i_1\in\mathbb I,~ \forall i_2...i_m\notin\mathbb I,\]
else, we say that $\A$ is irreducible.
\end{Definition}
\begin{Definition}\label{d1}\cite{A17}
A tensor $\A\in\Bbb R^{[m,n]}$ is called a $\Z$-tensor if its oﬀ-diagonal entries are non-positive. $\A$ is an $\M$-tensor if there
exists a nonnegative tensor $\B\in\Bbb R_+^{[m,n]}$ and a positive real number $\eta\geq\rho(\B)$ such that $\A=\eta\mathcal I-\B$. If $\eta>\rho(\B)$, then $\A$ is called a strong $\M$-tensor.\vspace{-0.3cm}
\end{Definition}
\begin{Definition}\label{d4}\cite{A8}
Let $A\in\Bbb R^{[2,n]}$ ($A$ is an $n$-dimensional square matrix) and $\B\in\Bbb R^{[m,n]}$. Then a product $\C=A\B\in\Bbb R^{[m,n]}$ is defined by
\begin{equation}\label{eq4}
c_{ji_2...i_m}=\sum_{j_2=1}^{n}a_{jj_2}b_{j_2i_2...i_m},
\end{equation}
which can be written as follows
\[\C_{(1)}=(A\B)_{(1)}=A\B_{(1)},\]
where $\C_{(1)}$ and $\B_{(1)}$ are the matrices obtained from $\C$ and $\B$ flattened along the first index, respectively.
\end{Definition}
\begin{Definition}\label{d2}\cite{A18}
Let $\A\in\Bbb R^{[m,n]}$. The majorization matrix of $\A$, denoted by $M(\A)$, is defined as a square matrix of size $n\times n$ with its entries
\[M(\A)_{ij}=a_{ij...j},~~i,j=1,2,...,n.\]
\end{Definition}
$~\hspace{0.5cm}$ If $M(\A)$ is a nonsingular matrix and
$\A=M(\A)\mathcal I_m$, then $M(\A)^{-1}$ is the order 2 left-inverse of $\A$, i.e., $M(\A)^{-1}\A=\I_m$, and then we call $\A$ a left-invertible tensor or left-nonsingular tensor.
\begin{Definition}\label{d3}\cite{A4}
Let $\A\in\Bbb R^{[m,n]}$. A pair $(\lambda,\textbf{x})\in\Bbb C\times(\Bbb C^n\backslash\{0\})$ is called an eigenvalue-eigenvector(or simply eigenpair) of $\A$ if they satisfy the equation
\begin{equation}\label{eq3}
\A\textbf{x}^{m-1}=\lambda\textbf{x}^{[m-1]},
\end{equation}
where $\textbf{x}^{[m-1]}=(x_1^{m-1},...,x_n^{m-1})^T$. We call $(\lambda,\textbf{x})$ an H-eigenpair if both $\lambda$ and $\textbf{x}$ are real.
\end{Definition}
$~\hspace{0.5cm}$Let $\rho(\A)=\max\{|\lambda| |\lambda\in\sigma(\A)\}$ be the spectral radius of $\A$, where $\sigma(\A)$ is the set of all eigenvalues of $\A$.
\begin{lem}\label{l1}\cite{A8}
If $\A$ is a strong $\M$-tensor, then $M(\A)$ is a nonsingular M-matrix.
\end{lem}
\begin{lem}\label{l2}\cite{A18}
If $M(\A)$ is an irreducible matrix, then $\A$ is irreducible.
\end{lem}
\begin{Definition}\label{d5}\cite{A8}
Let $\A,\mathcal E,\mathcal F\in\Bbb R^{[m,n]}$. $\A=\mathcal E-\mathcal F$ is said to be a splitting of $\A$ if $\mathcal E$ is a left-nonsingular; a regular splitting of $\A$ if $\mathcal E$ is left-nonsingular with
$M(\mathcal E)^{-1}\geq\textbf{O}$ and $\mathcal F\geq\mathcal O$; a weak regular splitting of $\A$ if $\mathcal E$ is left-nonsingular with
$M(\mathcal E)^{-1}\geq\textbf{O}$ and $M(\mathcal E)^{-1}\mathcal F\geq\mathcal O$; a convergent splitting if $\rho(M(\mathcal E)^{-1}\mathcal F)<1$.
\end{Definition}
\begin{lem}\label{l3}\cite{A22}
If $\A$ is a $\mathcal Z$-tensor, then the following conditions are equivalent
\begin{enumerate}
  \item $\A$ is a strong $\M$-tensor.
  \item $\A$ has a convergent (weak) regular splitting.
  \item All (weak) regular splittings of $\A$ are convergent.
  \item There exist a vector $\textbf{x}>\textbf{0}$ such that $\A\textbf{x}^{m-1}>\textbf{0}$.
\end{enumerate}
\end{lem}
\begin{lem}\label{l4}\cite{A3}
If $\A$ is a strong $\M$-tensor, then for every positive vector $\textbf{b}$, the multilinear system $\A\textbf{x}^{m-1}=\textbf{b}$ has a unique positive solution.
\end{lem}
\begin{lem}\label{l5}\cite{A20}
Suppose that $\A\in\Bbb R^{[m,n]}$. Let $\A=\mathcal E_1-\F_1=\mathcal E_2-\F_2$ be a weak regular splitting and a regular splitting, respectively, and $\F_2\leq\F_1,\F_2\neq\mathcal O$. One of the following statements holds.
\begin{enumerate}
\item $\rho(M(\mathcal E_2)^{-1}\F_2)\leq\rho(M(\mathcal E_1)^{-1}\F_1)<1$.
\item $\rho(M(\mathcal E_2)^{-1}\F_2)\geq\rho(M(\mathcal E_1)^{-1}\F_1)\geq1$.
\end{enumerate}
If $\F_2<\F_1,\F_2\neq\mathcal O$ and $\rho(M(\mathcal E_1)^{-1}\F_1)>1$, the first inequality in part $(2)$ is strict.
\end{lem}
\begin{lem}\label{l6}\cite{A22}
Let $\A$ be a strong $\M$-tensor, and $\A=\mathcal E_1-\F_1=\mathcal E_2-\F_2$ be two weak regular splitting with $M(\mathcal E_1)^{-1}\F_1\leq M(\mathcal E_2)^{-1}\F_2$. If the Perron vector $\textbf{x}$ of $M(\mathcal E_2)^{-1}\F_2$ satisfies $\A\textbf{x}^{m-1}\geq\textbf{0}$ then $\rho(M(\mathcal E_2)^{-1}\F_2)\leq\rho(M(\mathcal E_1)^{-1}\F_1)$.
\end{lem}
A general tensor splitting iterative method for solving \eqref{eq1} is
\begin{equation}\label{eq5}
\textbf{x}_{j+1}=[M(\mathcal E)^{-1}\mathcal F\textbf{x}_{j}^{m-1} +M(\mathcal E)^{-1}\textbf{b}]^{[\tfrac{1}{m-1}]},~
j=0,1,....
\end{equation}
$M(\mathcal E)^{-1}\mathcal F$ is called the iterative tensor of the splitting method \eqref{eq5}. Taking $\A=\D-\mathcal L-\F$, Liu et al. in \cite{A8}, considered $\E=\D,~\E=\D-\mathcal L$ and
$\E=\tfrac{1}{\tau}(\D-\tau\mathcal L)$, the Jacobian, the Gauss-Seidel, and the SOR iterative methods, respectively, where $\D={D}\I$ and $\mathcal L={L}\I$. ${D},~{L}$ are the positive diagonal matrix and the strictly lower triangle nonnegative matrix, respectively.
Without loss of generality, we always assume that $a_{ii...i}=1,~ i=1,2...,n$. Consider the splitting of $\A=\I-\mathcal L-\F$,
where $\mathcal L={L}\I$ and ${L}$ is the strictly lower triangle part of $M(\A)$.\\
\indent
Using iterative methods for solving \eqref{eq1} may have a poor convergence or even fail to converge. To overcome this problem, it is efficient to apply these methods which combine preconditioning techniques. These iterative methods usually involve some matrices that transform the iterative tensor $M(\mathcal E)^{-1}\mathcal F$ into a favorable tensor. The transformation matrices are called preconditioners. Li et al. in \cite{A10}, considered the preconditioner ${P}_{\boldsymbol{\alpha}}={I}+{S}_{\boldsymbol{\alpha}}$ for solving preconditioned multilinear system
\[{P}_{\boldsymbol{\alpha}}\A\textbf{x}^{m-1} ={P}_{\boldsymbol{\alpha}}\textbf{b},\]
with
\[{S}_{\boldsymbol{\alpha}}=\left[\begin{array}{*{20}c}
             0&-{\alpha}_1a_{12...2}&0&...&0\\
             0&0&-\alpha_2a_{23...3}&...&0\\
             \vdots & \vdots & \vdots &\ddots& \vdots \\
             0&0&0&...&\alpha_{n-1}a_{n-1,n...n}\\
             0&0&0&0&0
              \end{array}\right],\]
firstly proposed for $M$-matrix systems and the authors extended the results for solving tensor case. In \cite{A16}, Liu et al.
considered a new preconditioned SOR method for solving multilinear systems with preconditioner ${P}_{\boldsymbol{\beta}}={I}+{C}_{\boldsymbol{\beta}}$ where
\[{C}_{\boldsymbol{\beta}}=\left[\begin{array}{*{20}c}
                     0&0&0&...&0\\
                     -{\beta}_1a_{21...1}&0&0&...&0\\
                     \vdots & \vdots & \vdots &\ddots& \vdots \\
                     -{\beta}_{n-2}a_{(n-1)1...1}&0&0&...&0\\
                     -{\beta}_{n-1}a_{n1...1}&0&0&0&0
                     \end{array}\right].\]
\indent
Here we consider the preconditioner ${P}_{\boldsymbol{\alpha}\boldsymbol{\beta}}(s,k)=D+ {S}_{\boldsymbol{\alpha}}^{s}+{K}_{\boldsymbol{\beta}}^{k}$, where $1\leq s,k\leq n-1$, $D$ is the diagonal part of majorization of $\A$ (so herein $D=I$) and ${S}_{\boldsymbol{\alpha}}^{s}$, ${K}_{\boldsymbol{\beta}}^{k}$ are square matrices which all of their elements are zeros except the $s$th upper and the $k$th lower diagonals, i.e.,
\[{S}_{\boldsymbol{\alpha}}^{s}=\left[\begin{array}{*{20}c}
             0&...&0&-\alpha_1a_{1(1+s)...(1+s)}&0&...&0\\
             0&...&0&0&-\alpha_2a_{2(2+s)...(2+s)}&...&0\\
             \vdots&\vdots&\vdots&\vdots&\vdots&\ddots&0\\
             0&...&0&0&0&...&\alpha_{n-s}a_{n-s,n...n}\\
             0&...&0&0&0&...&0\\
             \vdots&\vdots&\vdots&\vdots&\vdots&\vdots&\vdots\\
             0&...&0&0&0&...&0
             \end{array}\right],\]
\[{K}_{\boldsymbol{\beta}}^{k}=\left[\begin{array}{*{20}c}
              0&0&...&0&...&0\\
              \vdots&\vdots&\vdots&\vdots&\vdots&\vdots\\
              0&0&...&0&...&0\\
             -{\beta}_{k+1}a_{(k+1)1...1}&0&...&0&...&0\\
             0&-{\beta}_{k+2}a_{(k+2)2...2}&...&0&...&0\\
             \vdots&\vdots&\ddots&\vdots&\vdots&\vdots\\
             0&0&-{\beta}_{n}a_{n(n-k+1)...(n-k+1)}&0&...&0
             \end{array}\right].\]
Applying ${P}_{\boldsymbol{\alpha}\boldsymbol{\beta}}(s,k)$ on the left side of Eq. \eqref{eq1},  we get a new preconditioned multi-linear system
\begin{equation}\label{eq6}
\A_{\boldsymbol{\alpha}\boldsymbol{\beta}}(s,k)\textbf{x}^{m-1} =\textbf{b}_{\boldsymbol{\alpha}\boldsymbol{\beta}}(s,k),
\end{equation}
with $\A_{\boldsymbol{\alpha}\boldsymbol{\beta}}(s,k)= {P}_{\boldsymbol{\alpha}\boldsymbol{\beta}}(s,k)\A$ and $\textbf{b}_{\boldsymbol{\alpha}\boldsymbol{\beta}}(s,k) ={P}_{\boldsymbol{\alpha}\boldsymbol{\beta}}(s,k)\textbf{b}.$
\begin{prop}\label{prop1}
Let $\A\in\Bbb{R}^{[m,n]}$ be a $\mathcal Z$-tensor. If $\A$ is a strong $\M$-tensor for any ${\beta}_j\in[0,1],~j=k+1,...,n$ and $\alpha_i\in[0,1],~i=1,...,n-s$, then $\A_{\boldsymbol{\alpha}\boldsymbol{\beta}}(s,k)$ is a strong $\mathcal M$-tensor.
\end{prop}
~\hspace{0.5cm}\textbf{Proof}.
Without loss of generality, we assume that $k=s=1$. Let $\A_{\boldsymbol{\alpha}\boldsymbol{\beta}}(s,k) ={P}_{\boldsymbol{\alpha}\boldsymbol{\beta}}(s,k)\A=(\hat{a}_{i_1i_2…i_m})$. Then for $1\leq i_2,…,i_m\leq n$, we have
\[\hat{a}_{ji_2…i_m}=\left\lbrace\begin{array}{l}
      {a}_{1i_2…i_m}-\alpha_1a_{12...2}{a}_{2i_2…i_m},~\hspace{7cm}j=1\\
{a}_{ji_2…i_m}-{\beta}_ja_{j(j-1)...(j-1)}{a}_{(j-1)i_2…i_m} -\alpha_ja_{j(j+1)...(j+1)}{a}_{(j+1)i_2…i_m},\hspace{0.7cm}2\leq j\leq n-1\\
{a}_{ni_2…i_m}-{\beta}_na_{n(n-1)...(n-1)}{a}_{(n-1)i_2…i_m},~\hspace{5cm}j=n.
                               \end{array}\right.\]
For $(j,i_2,…,i_m)\neq(j,j,...,j)$ and $\alpha_i,{\beta}_j\in[0,1]$, we have $\hat{a}_{ji_2…i_m}\leq0$, i.e., $\A_{\boldsymbol{\alpha}\boldsymbol{\beta}}(s,k)$ is a $\mathcal Z$-tensor. According to Lemma \ref{l3}, there exist a vector $\textbf{x}>\textbf{0}$ such that $\A\textbf{x}^{m-1}>\textbf{0}$. It follows from $P_{\boldsymbol{\alpha}}(s)>O$ that $P_{\boldsymbol{\alpha}}\A\textbf{x}^{m-1}>\textbf{0}$. Thus there exists a vector $\textbf{x}>\textbf{0}$ such that $\A_{\boldsymbol{\alpha}\boldsymbol{\beta}}(s,k)\textbf{x}^{m-1}>\textbf{0}$. Therefore, $\A_{\boldsymbol{\alpha}\boldsymbol{\beta}}(s,k)$ is a strong $\mathcal M$-tensor. $\hspace{1cm}\blacksquare$
\begin{prop}\label{prop2}
The preconditioned multi-linear system \eqref{eq6} has the same unique positive solution with multi-linear system \eqref{eq1}.
\end{prop}
~\hspace{0.5cm}\textbf{Proof}.
Because $\textbf{b}_{\boldsymbol{\alpha}\boldsymbol{\beta}}(s,k)\geq\textbf{b}>\textbf{0}$ for any ${\beta}_j\in[0,1],~j=k+1,...,n$ and $\alpha_i\in[0,1],~i=1,...,n-s$, by Lemma \ref{l4} and Proposition \ref{prop1}, it is obvious. $\hspace{2cm}\blacksquare$
\section{The preconditioned Jacobi, Gauss–Seidel and SOR type iteration schemes}
\subsection{\textbf{The preconditioned Jacobi type iterative scheme with the preconditioner $P_{\boldsymbol{\alpha}\boldsymbol{\beta}}(s,k)$}}
Let $\A=\I-\mathcal L-\F$. We consider the following five Jacobi type splittings:\vspace{0.2cm}\\
$\begin{array}{l}
\A_{\boldsymbol{\alpha}\boldsymbol{\beta}}(s,k) =P_{\boldsymbol{\alpha}\boldsymbol{\beta}}(s,k)\A\vspace{0.2cm}\\
\hspace{1.6cm}=P_{\boldsymbol{\alpha}\boldsymbol{\beta}}(s,k)\I -P_{\boldsymbol{\alpha}\boldsymbol{\beta}}(s,k)(\mathcal L+\mathcal F)=\mathcal E_1-\F_1.\vspace{0.2cm}\\
\A_{\boldsymbol{\alpha}\boldsymbol{\beta}}(s,k) =\I-(P_{\boldsymbol{\alpha}\boldsymbol{\beta}}(s,k)(\mathcal L+\mathcal F)-({S}_{\boldsymbol{\alpha}}^{s}+{K}_{\boldsymbol{\beta}}^{k})\I)=\mathcal E_2-\F_2.\vspace{0.2cm}\\
~~\A_{\boldsymbol{\alpha}}(s)=\I-(P_{\boldsymbol{\alpha}}(s)(\mathcal L+\mathcal F)-{S}_{\boldsymbol{\alpha}}^{s}\I)=\mathcal E_3-\F_3.\vspace{0.2cm}\\
~~\A_{\boldsymbol{\beta}}(k)=\I-(P_{\boldsymbol{\beta}}(\mathcal L+\mathcal F)-{K}_{\boldsymbol{\beta}}^{k}\I)=\mathcal E_4-\F_4.\vspace{0.3cm}\\
\end{array}$
\begin{rem}\label{rem1}
The splitting $\A_{\boldsymbol{\alpha}}(s)=\mathcal E_3-\F_3$, where $s=1$, is the same as the splitting in \cite{A10}.
\end{rem}
\begin{rem}\label{rem2}
When $s=k=1$, we denote ${K}_{\boldsymbol{\beta}}^{1}$ by ${K}_{\boldsymbol{\beta}}$ and ${S}_{\boldsymbol{\alpha}}^{1}$ by ${S}_{\boldsymbol{\alpha}}$. Thus we have the following Jacobi type splitting:\\
Denote ${K}_{\boldsymbol{\beta}}$ by ${K}$ and ${S}_{\boldsymbol{\alpha}}$ by ${S}$ when all ${\beta}_j=1,~j=2,3,...,n$ and ${\alpha}_i=1,~i=1,2,...,n-1$.\\
Let $\mathcal {L}={K}\I+\mathcal{{L'}}$ and  $\mathcal {F}={S}\I+\mathcal {{F'}}$, then\\
$\begin{array}{l}
\A_{\boldsymbol{\alpha}\boldsymbol{\beta}}(s,k) =(I-{S}_{\boldsymbol{\alpha}}K-{K}_{\boldsymbol{\beta}}{S})\I-[\mathcal L+\mathcal F-({S}_{\boldsymbol{\alpha}}+{K}_{\boldsymbol{\beta}})\I+ {S}_{\boldsymbol{\alpha}}(\mathcal{L'}+\mathcal {F})+{K}_{\boldsymbol{\beta}}(\mathcal{L}+\mathcal {F'})]=\mathcal E_5-\F_5.
\end{array}$
\end{rem}
\begin{prop}\label{prop3}
Let $\A\in\Bbb{R}^{[m,n]}$ be a strong $\M$-tensor for any $\boldsymbol{\beta}_j\in[0,1],~j=k+1,...,n$ and $\alpha_i\in[0,1],~i=1,...,n-s$, then $\A_{\boldsymbol{\alpha}\boldsymbol{\beta}}(s,k)=\mathcal E_1-\F_1=\mathcal E_2-\F_2$, $\A_{\boldsymbol{\alpha}}(s)=\mathcal E_3-\F_3$ and $\A_{\boldsymbol{\beta}}(k)=\mathcal E_4-\F_4$ are convergent. Moreover if
\begin{equation}\label{eq7}
\left\lbrace\begin{array}{l}
\displaystyle 0<\displaystyle \alpha_1\displaystyle a_{12...2}a_{21...1}<\displaystyle 1,\\
\displaystyle 0<\displaystyle \alpha_i\displaystyle a_{i(i+1)...(i+1)}\displaystyle a_{(i+1)i...i}+ \displaystyle{\beta}_i\displaystyle a_{i(i-1)...(i-1)a_{(i-1)i...i}}<\displaystyle 1,~i=2,...,n-1,\\
\displaystyle 0<\displaystyle \displaystyle{\beta}_n\displaystyle a_{n(n-1)...(n-1)a_{(n-1)n...n}}<\displaystyle 1,
\end{array}\right.
\end{equation}
then the tensor splitting $\A_{\boldsymbol{\alpha}\boldsymbol{\beta}}(s,k)=\mathcal E_5-\F_5$ is convergent.
\end{prop}
~\hspace{0.5cm}\textbf{Proof}.
Suppose $\A_{\boldsymbol{\alpha}\boldsymbol{\beta}}(s,k)=\mathcal E_1-\F_1.$ Since $\A=\I-\mathcal L-\F$ is a strong $\M$-tensor, $\rho(\mathcal L+\F)<1$. Thus $\rho(M(\mathcal E_1)^{-1}\F_1)=\rho(\mathcal L+\F)<1$. Hence $\A_{\boldsymbol{\alpha}\boldsymbol{\beta}}(s,k)=\mathcal E_1-\F_1$ is a convergent splitting.\\
Let $\A_{\boldsymbol{\alpha}\boldsymbol{\beta}}(s,k)=\mathcal E_2-\F_2$. We have $M(\mathcal E_1)^{-1}=I>O$ and since $\alpha_i,\boldsymbol{\beta}_j\in[0,1]$, it is easy to see that $\F_2\geq\mathcal O$. Thus $\A_{\boldsymbol{\alpha}\boldsymbol{\beta}}(s,k)=\mathcal E_2-\F_2$ is a regular splitting. By Proposition \ref{prop1}, $\A_{\boldsymbol{\alpha}\boldsymbol{\beta}}(s,k)$ is a strong $\M$-tensor and using Lemma \ref{l3}, $\A_{\boldsymbol{\alpha}\boldsymbol{\beta}}(s,k)=\mathcal E_2-\F_2$ is a convergent regular splitting.\\
When $\A_{\boldsymbol{\alpha}}(s)=\mathcal E_3-\F_3$ and $\A_{\boldsymbol{\beta}}(k)=\mathcal E_4-\F_4$, proof is similar to the proof of the case $\A_{\boldsymbol{\alpha}\boldsymbol{\beta}}(s,k)=\mathcal E_2-\F_2$.\\
Suppose that $\A_{\boldsymbol{\alpha}\boldsymbol{\beta}}(s,k)=\mathcal E_5-\F_5$, and Eq. \eqref{eq7}, holds. Thus $M(\mathcal E_5)^{-1}$ exists, and
\begin{equation}\label{eq8}
M(\mathcal E_5)_{ii}^{-1}=\left\lbrace\begin{array}{l}
\tfrac{\displaystyle1}{\displaystyle1-\displaystyle \alpha_1\displaystyle a_{12...2}a_{21...1}},\\
\tfrac{\displaystyle1}{\displaystyle1-\displaystyle \alpha_i\displaystyle a_{i(i+1)...(i+1)}a_{(i+1)i...i}- \displaystyle{\beta}_i\displaystyle a_{i(i-1)...(i-1)a_{(i-1)i...i}}},~i=2,...,n-1,\\
\tfrac{\displaystyle1}{\displaystyle1-\displaystyle \displaystyle{\beta}_n\displaystyle a_{n(n-1)...(n-1)a_{(n-1)n...n}}},
\end{array}\right.
\end{equation}
which implies that $M(\mathcal E_5)^{-1}>O$. It is not difficult to see that $\F_5=\mathcal E_5-\A_{\boldsymbol{\alpha}\boldsymbol{\beta}}(s,k)\geq\mathcal O$. Using Proposition \ref{prop1}, $\A_{\boldsymbol{\alpha}\boldsymbol{\beta}}(s,k)$ is a strong $\M$-tensor and from Lemma \ref{l3}, $\A_{\boldsymbol{\alpha}\boldsymbol{\beta}}(s,k)=\mathcal E_5-\F_5$ is a convergent regular splitting.$\hspace{2cm}\blacksquare$
\begin{prop}\label{prop4}
Let $\A$ be a strong $\M$-tensor and Eq. \eqref{eq7} holds. There exist $\textbf{x}_1,\textbf{x}_2\in\Bbb{R}_+^n$, such that
\begin{enumerate}
  \item $(M(\mathcal E_2)^{-1}\F_2)_{\boldsymbol{\alpha}\boldsymbol{\beta}} \textbf{x}_1^{m-1}\leq (M(\mathcal E_1)^{-1}\F_1)_ {\boldsymbol{\alpha}\boldsymbol{\beta}}\textbf{x}_1^{m-1}.$
 \item $\A_{\boldsymbol{\alpha}\boldsymbol{\beta}}(s,k) \textbf{x}_2^{m-1}\geq\textbf{0}$.
 \item $\rho((M(\mathcal E_5)^{-1}\F_5)_{\boldsymbol{\alpha}\boldsymbol{\beta}})\leq\rho((M(\mathcal E_2)^{-1}\F_2)_{\boldsymbol{\alpha}\boldsymbol{\beta}})\leq\rho((M(\mathcal E_1)^{-1}\F_1)_{\boldsymbol{\alpha}\boldsymbol{\beta}})$.
\end{enumerate}
\end{prop}
~\hspace{0.5cm}\textbf{Proof}.\\
1.~ Since $\A=\I-\mathcal L-\F$ is a strong $\M$-tensor, $\rho(M(\mathcal E_1)^{-1}(\mathcal L+\F))=\rho(\mathcal L+\F)<1$. Thus, for the nonnegative Jacobi iteration tensor $M(\mathcal E_1)^{-1}(\mathcal L+\F)=\mathcal L+\F$, there exists a nonnegative vector $\textbf{x}_1$ such that $M(\mathcal E_1)^{-1}(\mathcal L+\F)\textbf{x}_1^{m-1}=\rho(\mathcal L+\F)\textbf{x}_1^{[m-1]}$ by the Perron–Frobenius theorem. Thus we have\vspace{0.2cm}\\
$\begin{array}{l}
(M(\mathcal E_2)^{-1}\F_2)_{\boldsymbol{\alpha}\boldsymbol{\beta}}\textbf{x}_1^{m-1} =(P_{\boldsymbol{\alpha}\boldsymbol{\beta}}(s,k)(\mathcal L+\mathcal F)-({S}_{\boldsymbol{\alpha}}^{s}+{K}_{\boldsymbol{\beta}}^{k})\I)\textbf{x}_1^{m-1}\\
\hspace{3.45cm}=({I}+{S}_{\boldsymbol{\alpha}}^{s} +{K}_{\boldsymbol{\beta}}^{k})(\mathcal L+\mathcal F)\textbf{x}_1^{m-1}-({S}_{\boldsymbol{\alpha}}^{s} +{K}_{\boldsymbol{\beta}}^{k})\I\textbf{x}_1^{m-1}\\
\hspace{3.45cm}=(\mathcal L+\mathcal F)\textbf{x}_1^{m-1} +({K}_{\boldsymbol{\beta}}^{k}+{S} _{\boldsymbol{\alpha}}^{s})(\mathcal L+\mathcal F)\textbf{x}_1^{m-1} -({S}_{\boldsymbol{\alpha}}^{s} +{K}_{\boldsymbol{\beta}}^{k})\I\textbf{x}_1^{m-1}\\
\hspace{3.45cm}=(M(\mathcal E_1)^{-1}\F_1)_{\boldsymbol{\alpha}\boldsymbol{\beta}}\textbf{x}_1^{m-1} -({S}_{\boldsymbol{\alpha}}^{s}+{K}_{\boldsymbol{\beta}}^{k})(\I-(\mathcal L+\mathcal F))\textbf{x}_1^{m-1}\\
\hspace{3.45cm}=(M(\mathcal E_1)^{-1}\F_1)_{\boldsymbol{\alpha}\boldsymbol{\beta}}\textbf{x}_1^{m-1} -({S}_{\boldsymbol{\alpha}}^{s}+{K}_{\boldsymbol{\beta}}^{k})(1-\rho(M(\mathcal E_1)^{-1})\F_1)_{\boldsymbol{\alpha}\boldsymbol{\beta}}\textbf{x}_1^{[m-1]}.
 \end{array}$\vspace{0.2cm}\\
Thus\vspace{0.2cm}\\
$\begin{array}{l}
(M(\mathcal E_2)^{-1}\F_2)_{\boldsymbol{\alpha}\boldsymbol{\beta}} \textbf{x}_1^{m-1}-(M(\mathcal E_1)^{-1}\F_1)_ {\boldsymbol{\alpha}\boldsymbol{\beta}} \textbf{x}_1^{m-1} =-({S}_{\boldsymbol{\alpha}}^{s}+{K}_{\boldsymbol{\beta}}^{k})(1-\rho(M(\mathcal E_1)^{-1}\F_1))_{\boldsymbol{\alpha}\boldsymbol{\beta}} \textbf{x}_1^{[m-1]}\\
\hspace{5.9cm}\leq0,
\end{array}$\vspace{0.2cm}\\
due to ${S}_{\boldsymbol{\alpha}}^{s}+{K}_{\boldsymbol{\beta}}^{k}\geq O$ and $0<\rho(M(\mathcal E_1)^{-1}\F_1)<1$.\\
2.~ By Proposition \ref{prop3}, we know that $\A_{\boldsymbol{\alpha}\boldsymbol{\beta}}(s,k)=\mathcal E_5-\F_5$ is convergent, i.e., $0<\rho(M(\mathcal E_5)^{-1}\F_5)<1$ and thus, for the nonnegative Jacobi iteration tensor $M(\mathcal E_5)^{-1}\F_5$, there exists a nonnegative vector $\textbf{x}_2$ such that $M(\mathcal E_5)^{-1}\F_5\textbf{x}_2^{m-1}=\rho(M(\mathcal E_5)^{-1}\F_5)\textbf{x}_2^{[m-1]}$ by the Perron–Frobenius theorem. Thus we have\vspace{0.2cm}\\
$\begin{array}{l}
\A_{\boldsymbol{\alpha}\boldsymbol{\beta}}(s,k)\textbf{x}_2^{m-1}=\mathcal E_5\textbf{x}_2^{m-1} -\F_5\textbf{x}_2^{m-1}\\
\hspace{2.4cm}=\mathcal E_5\textbf{x}_2^{m-1} -M(\mathcal E_5)M(\mathcal E_5)^{-1}\F_5\textbf{x}_2^{m-1}\\
\hspace{2.4cm}=(I-{S}_{\boldsymbol{\alpha}}K-{K}_{\boldsymbol{\beta}}{S})\textbf{x}_2^{[m-1]} -\rho(M(\mathcal E_5)^{-1}\F_5)(I-{S}_{\boldsymbol{\alpha}}K-{K}_ {\boldsymbol{\beta}}{S})\I\textbf{x}_2^{m-1}\\
\hspace{2.4cm}=(I-{S}_{\boldsymbol{\alpha}}K-{K}_{\boldsymbol{\beta}}{S})\textbf{x}_2^{[m-1]} -\rho(M(\mathcal E_5)^{-1}\F_5)(I-{S}_{\boldsymbol{\alpha}}K-{K}_{\boldsymbol{\beta}}{S})\textbf{x}_2^{[m-1]}\\
\hspace{2.4cm}=(1-\rho(M(\mathcal E_5)^{-1}\F_5))(I-{S}_{\boldsymbol{\alpha}}K-{K} _{\boldsymbol{\beta}}{S})\textbf{x}_2^{[m-1]}\\
\hspace{2.4cm}\geq\textbf{0}.
\end{array}$\\
3.~Since $M(\mathcal E_2)_{\boldsymbol{\alpha}\boldsymbol{\beta}}^{-1}=I$ and $M(\mathcal E_5)_{\boldsymbol{\alpha}\boldsymbol{\beta}}^{-1}=(I-{S}_{\boldsymbol{\alpha}}K-{K}_{\boldsymbol{\beta}}{S})^{-1}$, thus $M(\mathcal E_5)_{\boldsymbol{\alpha}\boldsymbol{\beta}}^{-1}\geq M(\mathcal E_2)_{\boldsymbol{\alpha}\boldsymbol{\beta}}^{-1}$.\\
Let $(\rho((M(\mathcal E_5)^{-1}\F_5)_{\boldsymbol{\alpha}\boldsymbol{\beta}}),\textbf{x})$ be a Perron eigenpair of $(M(\mathcal E_5)^{-1}\F_5)_{\boldsymbol{\alpha}\boldsymbol{\beta}}$, then by part 2, we have $\A_{\boldsymbol{\alpha}\boldsymbol{\beta}}(s,k)\textbf{x}^{m-1}\geq\textbf{0}$ and by Lemma \ref{l6}, we have $\rho((M(\mathcal E_5)^{-1}\F_5)_{\boldsymbol{\alpha}\boldsymbol{\beta}})\leq\rho((M(\mathcal E_2)^{-1}\F_2)_{\boldsymbol{\alpha}\boldsymbol{\beta}})$. Now suppose that $\textbf{x}$ is a nonnegative Perron vector of $(M(\mathcal E_1)^{-1}\F_1)_{\boldsymbol{\alpha}\boldsymbol{\beta}}$, then by part 1, we have $(M(\mathcal E_2)^{-1}\F_2)_{\boldsymbol{\alpha}\boldsymbol{\beta}}\textbf{x}_1^{m-1}\leq (M(\mathcal E_1)^{-1}\F_1)_{\boldsymbol{\alpha}\boldsymbol{\beta}}\textbf{x}_1^{m-1}=\rho((M(\mathcal E_2)^{-1}\F_2)_{\boldsymbol{\alpha}\boldsymbol{\beta}})\textbf{x}_1^{[m-1]}.$ Since $(M(\mathcal E_2)^{-1}\F_2)_{\boldsymbol{\alpha}\boldsymbol{\beta}}\geq\mathcal O$, then $\rho((M(\mathcal E_2)^{-1}\F_2)_{\boldsymbol{\alpha}\boldsymbol{\beta}})\leq\rho((M(\mathcal E_1)^{-1}\F_1)_{\boldsymbol{\alpha}\boldsymbol{\beta}})$. This completes the proof.\hspace{2cm}$\blacksquare$
\begin{rem}\label{rem3}
It is easy to see that for every Perron vector of nonnegative Jacobi iteration tensor of convergence splitting method such as $\textbf{x}$, we have, $\A_{\boldsymbol{\alpha}\boldsymbol{\beta}}(s,k)\textbf{x}^{m-1}\geq\textbf{0}$.
\end{rem}
\begin{prop}\label{prop5} \cite{A22}
Let $\A\in\Bbb{R}^{[m,n]}$ be a strong $\M$-tensor. If Eq. \eqref{eq7} holds for any ${\beta}_{1,j},~{\beta}_{2,j}\in[0,1],~j=k+1,...,n$, $\alpha_{1,i},\alpha_{2,i}\in[0,1],~i=1,...,n-s$, $\boldsymbol{\alpha}'=(\alpha_{1,i}),~\boldsymbol{\alpha}''=(\alpha_{2,i}), ~\boldsymbol{\beta}'=({\beta}_{1,j}),~\boldsymbol{\beta}''=({\beta}_{2,j})$ and $\boldsymbol{\alpha}'\geq\boldsymbol{\alpha}'',~\boldsymbol{\beta}'\geq\boldsymbol{\beta}''$, then we have
\begin{enumerate}
  \item $\rho((M(\mathcal E_1)^{-1}\F_1)_{\boldsymbol{\alpha}'\boldsymbol{\beta}'})\leq \rho((M(\mathcal E_1)^{-1}\F_1)_{\boldsymbol{\alpha}''\boldsymbol{\beta}''})$.
  \item $\rho((M(\mathcal E_2)^{-1}\F_2)_{\boldsymbol{\alpha}'\boldsymbol{\beta}'})\leq \rho((M(\mathcal E_2)^{-1}\F_2)_{\boldsymbol{\alpha}''\boldsymbol{\beta}''})$.
  \item $\rho((M(\mathcal E_3)^{-1}\F_3)_{\boldsymbol{\alpha}'})\leq \rho((M(\mathcal E_3)^{-1}\F_3)_{\boldsymbol{\alpha}''})$.
  \item $\rho((M(\mathcal E_4)^{-1}\F_4)_{\boldsymbol{\beta}'})\leq \rho((M(\mathcal E_4)^{-1}\F_4)_{\boldsymbol{\beta}''})$.
  \item $\rho((M(\mathcal E_5)^{-1}\F_5)_{\boldsymbol{\alpha}'\boldsymbol{\beta}'})\leq \rho((M(\mathcal E_5)^{-1}\F_5)_{\boldsymbol{\alpha}''\boldsymbol{\beta}''})$.
\end{enumerate}
\end{prop}
\subsection{\textbf{Gauss-Seidel type iterative schemes with the tridiagonal preconditioner $P_{\boldsymbol{\alpha}\boldsymbol{\beta}}(s,k)$}}
We consider the following four Gauss-Seidel type splittings:\vspace{0.2cm}\\
$\begin{array}{l}
\A_{\boldsymbol{\alpha}\boldsymbol{\beta}}(s,k)= P_{\boldsymbol{\alpha}\boldsymbol{\beta}}(s,k)\A\vspace{0.2cm}\\
\hspace{1.6cm}=P_{\boldsymbol{\alpha}\boldsymbol{\beta}}(s,k)(\I-\mathcal L)-P_{\boldsymbol{\alpha}\boldsymbol{\beta}}(s,k)\mathcal F=\M_1-\N_1.\vspace{0.2cm}\\
\A_{\boldsymbol{\alpha}\boldsymbol{\beta}}(s,k)=(\I-\mathcal L+K_{\boldsymbol{\beta}}^k\I-K_{\boldsymbol{\beta}}^k\mathcal L-\D_{\boldsymbol{\alpha}}-\mathcal L_{\boldsymbol{\alpha}}-\D_{\boldsymbol{\beta}}-\mathcal L_{\boldsymbol{\beta}}) -(\F-S_{\boldsymbol{\alpha}}^s\I+S_{\boldsymbol{\alpha}}^s\F+ \F_{\boldsymbol{\alpha}}+\F_{\boldsymbol{\beta}})=\M_2-\N_2.\vspace{0.2cm}\\
\A_{\boldsymbol{\alpha}}(s)=(\I-\mathcal L-\D_{\boldsymbol{\alpha}}-\mathcal L_{\boldsymbol{\alpha}}) -(\F-S_{\boldsymbol{\alpha}}^s\I+S_{\boldsymbol{\alpha}}^s\F +\F_{\boldsymbol{\alpha}})=\M_3-\N_3.\vspace{0.2cm}\\
\A_{\boldsymbol{\beta}}(k)=((I+K_{\boldsymbol{\beta}}^k)(\I-\mathcal L)-\D_{\boldsymbol{\beta}}-\mathcal L_{\boldsymbol{\beta}})-(\F+\F_{\boldsymbol{\beta}})=\M_4-\N_4.\vspace{0.3cm}\\
\end{array}$\\
Where $\D_{\boldsymbol{\alpha}}=D_{\boldsymbol{\alpha}}\I,~\mathcal L_{\boldsymbol{\alpha}}=L_{\boldsymbol{\alpha}}\I, ~\D_{\boldsymbol{\beta}}=D_{\boldsymbol{\beta}}\I,~\mathcal L_{\boldsymbol{\beta}}=L_{\boldsymbol{\beta}}\I$, and $D_{\boldsymbol{\alpha}},~D_{\boldsymbol{\beta}}, ~L_{\boldsymbol{\alpha}},~L_{\boldsymbol{\beta}}$ are the diagonal parts and the strictly lower triangle parts of $M(S_{\boldsymbol{\alpha}}^s\mathcal L)$ and $M(K_{\boldsymbol{\beta}}^k\F)$, respectively, i.e.
\[S_{\boldsymbol{\alpha}}^s\mathcal L=\D_{\boldsymbol{\alpha}}+\mathcal L_{\boldsymbol{\alpha}}+\F_{\boldsymbol{\alpha}},~ K_{\boldsymbol{\beta}}^k\F=\D_{\boldsymbol{\beta}}+\mathcal L_{\boldsymbol{\beta}}+\F_{\boldsymbol{\beta}}.\]
\begin{rem}\label{rem4}
The splitting $\A_{\boldsymbol{\alpha}}(s)=\M_3-\N_3$, where $s=1$, is the same as the splitting in \cite{A10}.
\end{rem}
\begin{rem}\label{rem5}
If $k=l=1$, similar to Remark \ref{rem2}, we have
\[
\A_{\boldsymbol{\alpha}\boldsymbol{\beta}}(s,k)= ((I+{K}_{\boldsymbol{\beta}})(\I-\mathcal L)-S_{\boldsymbol{\alpha}}
\mathcal L- K_{\boldsymbol{\beta}}S\I)-((I+S_{\boldsymbol{\alpha}})\F -S_{\boldsymbol{\alpha}}\I+K_{\boldsymbol{\beta}}\F') =\mathcal M_5 -\mathcal N_5.
\]
\end{rem}
\begin{prop}\label{prop6}
Let $\A\in\Bbb{R}^{[m,n]}$ be a strong $\M$-tensor for any $\boldsymbol{\beta}_j\in[0,1],~j=k+1,...,n$ and $\alpha_i\in[0,1],~i=1,...,n-s$, then $\A_{\boldsymbol{\alpha}\boldsymbol{\beta}}(s,k)=\M_1-\N_1$ is convergent.\\
When $k<s$ if
\begin{equation}\label{eq9}
\left\lbrace\begin{array}{l}
0<\displaystyle\alpha_i\displaystyle a_{i(n-i)...(n-i)}\displaystyle a_{(n-i)i...i}<1,\hspace{4.5cm}~i=1,2,...,k,\\
\displaystyle 0<\displaystyle\alpha_i\displaystyle a_{i(n-i)...(n-i)}\displaystyle a_{(n-i)i...i}+ \displaystyle{\beta}_i\displaystyle a_{i(i-k)...(i-k)}\displaystyle a_{(i-k)i...i}<\displaystyle 1,\hspace{0.5cm}i=k+1,...,s,\\
0<\displaystyle{\beta}_i\displaystyle a_{i(i-k)...(i-k)}\displaystyle a_{(i-k)i...i}<1,\hspace{4.5cm}~i=s+1,...,n.
\end{array}\right.
\end{equation}
When $k>s$ if
\begin{equation}\label{eq10}
\left\lbrace\begin{array}{l}
0<\displaystyle\alpha_i\displaystyle a_{i(n-i)...(n-i)}\displaystyle a_{(n-i)i...i}<1,\hspace{4.5cm}~i=1,2,...,s,\\
\displaystyle 0<\displaystyle\alpha_i\displaystyle a_{i(n-i)...(n-i)}\displaystyle a_{(n-i)i...i}+ \displaystyle{\beta}_i\displaystyle a_{i(i-k)...(i-k)}\displaystyle a_{(i-k)i...i}<\displaystyle 1,\hspace{0.5cm}i=s+1,...,k,\\
0<\displaystyle{\beta}_i\displaystyle a_{i(i-k)...(i-k)}\displaystyle a_{(i-k)i...i}<1,\hspace{4.5cm}~i=k+1,...,n.
\end{array}\right.
\end{equation}
And when $k=s$ if
\begin{equation}\label{eq11}
\left\lbrace\begin{array}{l}
0<\displaystyle\alpha_i\displaystyle a_{i(i+k)...(i+k)}\displaystyle a_{(i+k)i...i}<1,\hspace{4.5cm}~i=1,2,...,k,\\
\displaystyle 0<\displaystyle\alpha_i\displaystyle a_{i(i+k)...(i+k)}\displaystyle a_{(i+k)i...i}+ \displaystyle{\beta}_i\displaystyle a_{i(i-k)...(i-k)}\displaystyle a_{(i-k)i...i}<\displaystyle 1,\hspace{0.5cm}i=k+1,...,n-k,\\
0<\displaystyle{\beta}_i\displaystyle a_{i(i-k)...(i-k)}\displaystyle a_{(i-k)i...i}<1,\hspace{4.5cm}~i=n-k+1,...,n.
\end{array}\right.
\end{equation}
Then the tensor splitting $\A_{\boldsymbol{\alpha}\boldsymbol{\beta}}(s,k)=\M_2-\N_2$, is convergent.
If \[0<\displaystyle\alpha_i\displaystyle a_{i(n-i)...(n-i)}\displaystyle a_{(n-i)i...i}<1,~i=1,2,...,s,\]
$\A_{\boldsymbol{\alpha}}(s)=\M_3-\N_3$ is convergent. Finally, if
\[0<\displaystyle{\beta}_i\displaystyle a_{i(i-k)...(i-k)}\displaystyle a_{(i-k)i...i}<1,~i=k+1,...,n,\]
then $\A_{\boldsymbol{\beta}}(k)=\M_4-\N_4$ is convergent.
\end{prop}
~\hspace{0.5cm}\textbf{Proof}. Let $\A_{\boldsymbol{\alpha}\boldsymbol{\beta}}(s,k)=\M_1-\N_1$. Due to Proposition \ref{prop1}, $\A_{\boldsymbol{\alpha}\boldsymbol{\beta}}(s,k)$ is a strong $\M$-tensor, and $\N_1\geq\mathcal O$. Since
\[M(\M_1)^{-1}\N_1=(I-L)^{-1}P_{\boldsymbol{\alpha}\boldsymbol{\beta}}(s,k)^{-1}P_{\boldsymbol{\alpha}\boldsymbol{\beta}}(s,k)\N_1\geq\mathcal O,\]
$\A_{\boldsymbol{\alpha}\boldsymbol{\beta}}(s,k)=\M_1-\N_1$ is a weak regular splitting and using Lemma \ref{l3}, is convergent.\\
Suppose that $\A_{\boldsymbol{\alpha}\boldsymbol{\beta}}(s,k)=\M_2-\N_2$ and $k=s.$ Since $\M_2=\I-\mathcal L+K_{\boldsymbol{\beta}}^k\I-K_{\boldsymbol{\beta}}^k\mathcal L-\D_{\boldsymbol{\alpha}}-\mathcal L_{\boldsymbol{\alpha}}-\D_{\boldsymbol{\beta}}-\mathcal L_{\boldsymbol{\beta}}$, then $M(\M_2)=I-D_{\boldsymbol{\alpha}}-D_{\boldsymbol{\beta}}-L+K_{\boldsymbol{\beta}}^k-K_{\boldsymbol{\beta}}^kL-L_{\boldsymbol{\alpha}}-L_{\boldsymbol{\beta}}$.
Notice that $D_{\boldsymbol{\alpha}}$ and $D_{\boldsymbol{\beta}}$ are diagonal part of $M(S_{\boldsymbol{\alpha}}^s\mathcal L)$ and $M(K_{\boldsymbol{\beta}}^k\F)$, respectively. It is not difficult to see that
\begin{equation}\label{eq12}
(I-D_{\boldsymbol{\alpha}}-D_{\boldsymbol{\beta}})_{ii} =\left\lbrace\begin{array}{l}
{1-\displaystyle{\alpha}_i\displaystyle a_{i(i+k)...(i+k)}\displaystyle a_{(i+k)i...i}},\hspace{4.5cm}~i=1,2,...,k,\\
{1-\displaystyle\alpha_i\displaystyle a_{i(i+k)...(i+k)}\displaystyle a_{(i+k)i...i}- \displaystyle{\beta}_i\displaystyle a_{i(i-k)...(i-k)}\displaystyle a_{(i-k)i...i}},\hspace{0.5cm}i=k+1,...,n-k,\\
{1-\displaystyle{\beta}_i\displaystyle a_{i(i-k)...(i-k)}\displaystyle a_{(i-k)i...i}},\hspace{4.5cm}~i=n-k+1,...,n.
\end{array}\right.
\end{equation}
Since Eq. \eqref{eq11} holds, $(I-D_{\boldsymbol{\alpha}}-D_{\boldsymbol{\beta}})^{-1}$ exists, and $(I-(D_{\boldsymbol{\alpha}}+D_{\boldsymbol{\beta}}))^{-1}=I+(D_{\boldsymbol{\alpha}}+D_{\boldsymbol{\beta}})+... +(D_{\boldsymbol{\alpha}}+D_{\boldsymbol{\beta}})^{n-1}\geq I$. Denote $H:=L+L_{\boldsymbol{\alpha}}+L_{\boldsymbol{\beta}}-K_{\boldsymbol{\beta}}^k+K_{\boldsymbol{\beta}}^kL.$ $H$ is a lower triangular matrix, to prove $H\geq O$ it is sufﬁcient to show $(L-K_{\boldsymbol{\beta}}^k)_{i+1,i}\geq0$ for any $i=1,...,n-1$. Actually,
\[(L-K_{\boldsymbol{\beta}}^k)_{i+1,i}=-a_{i+1,i...i}-(-\boldsymbol{\beta}_{i+1}a_{i+1,i...i})= a_{i+1,i...i}(\boldsymbol{\beta}_{i+1}-1)\geq0.\]
By the Neumann’s series \cite{A23}, we have\vspace{0.2cm}\\
$\begin{array}{l}
M(\M_2)^{-1}=[(I-D_{\boldsymbol{\alpha}}-D_{\boldsymbol{\beta}})-H]^{-1}\\
\hspace{1.7cm}=[I-(I-(D_{\boldsymbol{\alpha}}+D_{\boldsymbol{\beta}}))^{-1}H]^{-1} (I-(D_{\boldsymbol{\alpha}}+D_{\boldsymbol{\beta}}))^{-1}\\
\hspace{1.7cm}=\{I+(I-(D_{\boldsymbol{\alpha}}+D_{\boldsymbol{\beta}}))^{-1}H +[(I-(D_{\boldsymbol{\alpha}}+D_{\boldsymbol{\beta}})^{-1}H]^2+...+\\
\hspace{3cm}+[(I-(D_{\boldsymbol{\alpha}}+D_{\boldsymbol{\beta}}))^{-1}H]^{n-1}\} (I-(D_{\boldsymbol{\alpha}}+D_{\boldsymbol{\beta}}))^{-1}\\
\hspace{1.7cm}\geq O.
\end{array}$\vspace{0.2cm}\\
Since $\N_2\geq\mathcal O$ (by the same discussion in proofing how $H\geq O$), $\A_{\boldsymbol{\alpha}\boldsymbol{\beta}}(s,k)=\M_2-\N_2$ is a weak regular splitting and using Lemma \ref{l3}, is convergent. For cases $k<s$ and $k>s$, similar discussion can be used for obtain desired results.\\
When $\A_{\boldsymbol{\alpha}}(s)=\M_3-\N_3$ and $\A_{\boldsymbol{\beta}}(k)=\M_4-\N_4$, proof is similar to the proof of the case $\A_{\boldsymbol{\alpha}\boldsymbol{\beta}}(s,k)=\M_2-\N_2$. $\hspace{2cm}\blacksquare$
\begin{prop}\label{prop7}
Let $\A$ be a strong $\M$-tensor and Eqs. \eqref{eq9}-\eqref{eq11} hold. There exists $\textbf{x}\in\Bbb{R}_+^n$, such that
\begin{enumerate}
 \item $\A_{\boldsymbol{\alpha}\boldsymbol{\beta}}(s,k) \textbf{x}^{m-1}\geq\textbf{0}$.
 \item $\rho((M(\M_2)^{-1}\N_2)_{\boldsymbol{\alpha} \boldsymbol{\beta}})\leq\rho((M(\M_3)^{-1}\N_3) _{\boldsymbol{\alpha}\boldsymbol{\beta}})\leq \rho((M(\M_1)^{-1}\N_1)_{\boldsymbol{\alpha}\boldsymbol{\beta}})<1$.
  \item $\rho((M(\M_2)^{-1}\N_2)_{\boldsymbol{\alpha} \boldsymbol{\beta}})\leq\rho((M(\M_4)^{-1}\N_4) _{\boldsymbol{\alpha}\boldsymbol{\beta}})\leq \rho((M(\M_1)^{-1}\N_1)_{\boldsymbol{\alpha}\boldsymbol{\beta}})<1$.
\end{enumerate}
\end{prop}
~\hspace{0.5cm}\textbf{Proof}.\\
1.~$\A_{\boldsymbol{\alpha}\boldsymbol{\beta}}(s,k)$ is a strong $\M$-tensor by Proposition \ref{prop1}. Using Lemma \ref{l3}, there exists $\textbf{x}\in\Bbb{R}_+^n$ such that $\A_{\boldsymbol{\alpha}\boldsymbol{\beta}}(s,k)\textbf{x}^{m-1}\geq\textbf{0}$.\\
2.~ From Proposition \ref{prop6}, $\A_{\boldsymbol{\alpha}\boldsymbol{\beta}}(s,k)=\M_2-\N_2$ and $\A_{\boldsymbol{\alpha}\boldsymbol{\beta}}(s,k)=\M_3-\N_3$ are two weak regular splitting. Denote $H':=L+L_{\boldsymbol{\alpha}}\geq O$. By Neumann’s series \cite{A23}, we have \vspace{0.2cm}\\
$\begin{array}{l}
M(\M_2)^{-1}=[(I-D_{\boldsymbol{\alpha}}-D_{\boldsymbol{\beta}})-H]^{-1}\\
\hspace{1.7cm}=[I-(I-(D_{\boldsymbol{\alpha}}+D_{\boldsymbol{\beta}}))^{-1}H]^{-1} (I-(D_{\boldsymbol{\alpha}}+D_{\boldsymbol{\beta}}))^{-1}\\
\hspace{1.7cm}=\{I+(I-(D_{\boldsymbol{\alpha}}+D_{\boldsymbol{\beta}}))^{-1}H +[(I-(D_{\boldsymbol{\alpha}}+D_{\boldsymbol{\beta}})^{-1}H]^2+...+\\
\hspace{3cm}+[(I-(D_{\boldsymbol{\alpha}}+D_{\boldsymbol{\beta}}))^{-1}H]^{n-1}\} (I-(D_{\boldsymbol{\alpha}}+D_{\boldsymbol{\beta}}))^{-1}\\
\hspace{1.7cm}\geq\{I+(I-D_{\boldsymbol{\alpha}})^{-1}H' +[(I-D_{\boldsymbol{\alpha}})^{-1}H']^2+...+
[(I-D_{\boldsymbol{\alpha}})^{-1}H']^{n-1}\} (I-D_{\boldsymbol{\alpha}})^{-1}\\
\hspace{1.7cm}=[(I-D_{\boldsymbol{\alpha}})-H']^{-1}\\
\hspace{1.7cm}=M(\M_1)^{-1}.
\end{array}$\\
By Proposition \ref{prop6}, we know that $\A_{\boldsymbol{\alpha}\boldsymbol{\beta}}(s,k)=\M_2-\N_2$ is convergent, i.e., $0<\rho(M(\M_2)^{-1}\N_2)<1$ and thus, for the nonnegative Gauss-Seidel iteration tensor $M(\M_2)^{-1}\N_2$, there exists a nonnegative vector $\textbf{x}$ such that $M(\M_2)^{-1}\N_2\textbf{x}^{m-1}=\rho(M(\M_2)^{-1}\N_2)\textbf{x}^{[m-1]}$ by the Perron–Frobenius theorem. Using Proposition \ref{prop4}, we have $\A_{\boldsymbol{\alpha}\boldsymbol{\beta}}(s,k)\textbf{x}^{m-1}\geq0$. By Lemma \ref{l6}, we have $\rho((M(\M_2)^{-1}\N_2)_{\boldsymbol{\alpha}\boldsymbol{\beta}})\leq\rho((M(\M_3)^{-1}\N_3) _{\boldsymbol{\alpha}\boldsymbol{\beta}})$. Similar discussion give us $\rho((M(\M_3)^{-1}\N_3) _{\boldsymbol{\alpha}\boldsymbol{\beta}})\leq\rho((M(\M_1)^{-1}\N_1)_{\boldsymbol{\alpha}\boldsymbol{\beta}})$. Since $\A_{\boldsymbol{\alpha}\boldsymbol{\beta}}(s,k)=\M_1-\N_1$ is convergent by Proposition \ref{prop6}, thus $\rho((M(\M_1)^{-1}\N_1)_{\boldsymbol{\alpha}\boldsymbol{\beta}})<1$.\\
3.~The proof of part 3 is similar to the proof of previous part.$\hspace{2cm}\blacksquare$
\begin{prop}\label{prop8}
Let $\A\in\Bbb{R}^{[m,n]}$ be a strong $\M$-tensor. If Eqs. \eqref{eq9}-\eqref{eq11} hold for any\\ ${\beta}_{1,j},~{\beta}_{2,j}\in[0,1],~j=k+1,...,n$, $\alpha_{1,i},\alpha_{2,i}\in[0,1],~i=1,...,n-s$, $\boldsymbol{\alpha}'=(\alpha_{1,i}),~\boldsymbol{\alpha}''=(\alpha_{2,i}),\\ ~\boldsymbol{\beta}'=({\beta}_{1,j}), ~\boldsymbol{\beta}''=({{\beta}}_{2,j})$ and $\boldsymbol{\alpha}'\geq\boldsymbol{\alpha}'', ~\boldsymbol{\beta}'\geq\boldsymbol{\beta}''$, then
\begin{enumerate}
  \item $\rho((M(\M)^{-1}\N_1)_{\boldsymbol{\alpha}'\boldsymbol{\beta}'})\leq \rho((M(\M_1)^{-1}\N_1)_{\boldsymbol{\alpha}''\boldsymbol{\beta}''})$.
  \item $\rho((M(\M_2)^{-1}\N_2)_{\boldsymbol{\alpha}'\boldsymbol{\beta}'})\leq \rho((M(\M_2)^{-1}\N_2)_{\boldsymbol{\alpha}''\boldsymbol{\beta}''})$.
  \item $\rho((M(\M_3)^{-1}\N_3)_{\boldsymbol{\alpha}'})\leq \rho((M(\M_3)^{-1}\N_3)_{\boldsymbol{\alpha}''})$.
  \item $\rho((M(\M_4)^{-1}\N_4)_{\boldsymbol{\beta}'})\leq \rho((M(\M_4)^{-1}\N_4)_{\boldsymbol{\beta}''})$.
\end{enumerate}
\end{prop}
\subsection{\textbf{The preconditioned SOR type method}}
In \cite{A8}, the SOR type method for solving Eq. \eqref{eq1} is given by taking $\mathcal E=\tfrac{1}{\displaystyle\omega}(\I-\displaystyle \omega\mathcal L)$ and
\[\textbf{x}_{j+1}=(M(\I-\displaystyle\omega\mathcal L)^{-1} ((1-\displaystyle\omega)\I+\displaystyle\omega\F)\textbf{x}_j^{m-1} +\displaystyle\omega M(\I-\omega\mathcal L)^{-1}\textbf{b})^{[\tfrac{1}{m-1}]}.\]
In this paper, we consider the following preconditioned SOR type method: ~\hspace{0.2cm}\\
$\begin{array}{l}
~\hspace{1cm}\textbf{x}_{j+1}=(\mathcal H_{\boldsymbol{\alpha}\boldsymbol{\beta}}(\displaystyle\omega) \textbf{x}_j^{m-1} +\textbf{h}_{\boldsymbol{\alpha}\boldsymbol{\beta}} ({\displaystyle\omega}))^{[\tfrac{1}{m-1}]},
\end{array}$~\hspace{0.2cm}\\
where~\hspace{0.2cm}\\
$\begin{array}{l}
~\hspace{1cm}\mathcal H_{\boldsymbol{\alpha}\boldsymbol{\beta}}(\displaystyle\omega)=M(\mathcal E_{\boldsymbol{\alpha}\boldsymbol{\beta}} (\displaystyle\omega))^{-1} \F_{\boldsymbol{\alpha}\boldsymbol{\beta}}(\displaystyle\omega), ~\textbf{h}_{\boldsymbol{\alpha}\boldsymbol{\beta}}(\displaystyle\omega)=M(\mathcal E_{\boldsymbol{\alpha}\boldsymbol{\beta}} (\displaystyle\omega))^{-1}\textbf{b}_{\boldsymbol{\alpha}\boldsymbol{\beta}}(s,k),\\
~\hspace{1cm}\mathcal E_{\boldsymbol{\alpha}\boldsymbol{\beta}}({\displaystyle\omega})= \tfrac{1}{\displaystyle\omega}(\D_{\boldsymbol{\alpha}\boldsymbol{\beta}} -\displaystyle\omega\mathcal L_{\boldsymbol{\alpha}\boldsymbol{\beta}}), ~\F_{\boldsymbol{\alpha}\boldsymbol{\beta}}({\displaystyle\omega})= \tfrac{1}{\displaystyle\omega} ((1-\displaystyle\omega)\D_{\boldsymbol{\alpha}\boldsymbol{\beta}}+ \displaystyle\omega\F_{\boldsymbol{\alpha}\boldsymbol{\beta}}),
\end{array}$\\
and\vspace{0.2cm}\\
$\begin{array}{l}
~\hspace{1cm}\D_{\boldsymbol{\alpha}\boldsymbol{\beta}}= \I-\D_{\boldsymbol{\alpha}}-\D_{\boldsymbol{\beta}},\vspace{0.2cm}\\
~\hspace{1cm}\mathcal L_{\boldsymbol{\alpha}\boldsymbol{\beta}}=\mathcal L -K_{\boldsymbol{\beta}}^k\I+K_{\boldsymbol{\beta}}^k\mathcal L+\mathcal L_{\boldsymbol{\alpha}}+\mathcal L_{\boldsymbol{\beta}},\vspace{0.2cm}\\
~\hspace{1cm}\F_{\boldsymbol{\alpha}\boldsymbol{\beta}} =\F-S_{\boldsymbol{\alpha}}^s\I+S_{\boldsymbol{\alpha}}^s\F +\F_{\boldsymbol{\alpha}}+\F_{\boldsymbol{\beta}}. \vspace{0.2cm}
\end{array}$
\begin{rem}\label{rem1}
When $s=1$ and $k=0$, the new preconditioned SOR method is similar to the preconditioned SOR method which is proposed in \cite{A10}.
\end{rem}
\begin{prop}\label{prop9}
Let $\A\in\Bbb{R}^{[m,n]}$ be a strong $\M$-tensor. If $\A=\I-\mathcal L-\F$ and $0<\omega_1<\omega_2\leq1$, then $\rho(\mathcal H_{\boldsymbol{\alpha}\boldsymbol{\beta}}(\displaystyle\omega_2))\leq\rho(\mathcal H_{\boldsymbol{\alpha}\boldsymbol{\beta}}(\displaystyle\omega_1))<1$.
\end{prop}
\begin{prop}\label{prop10}
Let $\A\in\Bbb{R}^{[m,n]}$ be a strong $\M$-tensor. For any $\omega\in(0,1],~\rho(\Theta_{\boldsymbol{\alpha}\boldsymbol{\beta}})\leq\rho(\mathcal H_{\boldsymbol{\alpha}\boldsymbol{\beta}}(\displaystyle\omega))$, where $\Theta_{\boldsymbol{\alpha}\boldsymbol{\beta}}$ is the iteration tensor of the preconditioned Gauss-Seidel type methods.
\end{prop}
\begin{prop}\label{prop11}\cite{A16}
Let $\A\in\Bbb{R}^{[m,n]}$ be a strong $\M$-tensor with $a_{i(i+1)...(i+1)}a_{(i+1)i...i}>0,~i=1,2,...,n-1$ and $0<a_{i1...1}a_{1i...i}<1,~i=2,3,...,n$. Then $\mathcal H(\omega)$ is nonnegative and irreducible for $0<\omega<1$.
\end{prop}
\begin{prop}\label{prop12}\cite{A16}
Let $\A\in\Bbb{R}^{[m,n]}$ be a strong $\M$-tensor and ${\alpha}_i,{\beta}_j\in[0,1],~i=1,2,...,n-1$. Then
If $0<\omega\leq1$, $a_{i(i+1)...(i+1)}a_{(i+1)i...i}>0,~i=1,2,...,n-1$ and $0<a_{i1...1}a_{1i...i}<1,~i=2,3,...,n$, we have
\[\rho(\mathcal H_{\boldsymbol{\alpha}\boldsymbol{\beta}}(\omega))\leq\rho(\mathcal H(\omega))<1.\]
\end{prop}
\begin{prop}\label{prop13}\cite{A16}
Let $\A\in\Bbb{R}^{[m,n]}$ be a strong $\M$-tensor. If the conditions of the Proposition \ref{prop11} hold and
${\beta}_{1,j},~{\beta}_{2,j}\in[0,1],~j=k+1,...,n$, ${\alpha}_{1,i},{\alpha}_{2,i}\in[0,1],~i=1,...,n-s$, $\boldsymbol{\alpha}'=(\alpha_{1,i}),~\boldsymbol{\alpha}''=(\alpha_{2,i}), ~\boldsymbol{\beta}'=({\beta}_{1,j}),~\boldsymbol{\beta}''=({{\beta}}_{2,j})$ and $\boldsymbol{\alpha}'\geq\boldsymbol{\alpha}'', ~\boldsymbol{\beta}'\geq\boldsymbol{\beta}''$, then
\[\rho(\mathcal H_{\boldsymbol{\alpha}'\boldsymbol{\beta}'}(\omega))\leq\rho(\mathcal H_{\boldsymbol{\alpha}''\boldsymbol{\beta}''}(\omega))<1.\]
\end{prop}
\begin{prop}\label{prop14}\cite{A16}
Let $\A\in\Bbb{R}^{[m,n]}$ be a strong $\M$-tensor. If $0<\omega_1<\omega_2\leq1$, then
\[\rho(M(\I-\displaystyle\omega_2\mathcal L)^{-1} (\omega_2\F+(1-\displaystyle\omega_2)\I))\leq \rho(M(\I-\displaystyle\omega_1\mathcal L)^{-1} (\omega_1\F+(1-\displaystyle\omega_1)\I))<1.\]
\end{prop}
\begin{prop}\label{prop15}\cite{A16}
Let $\A\in\Bbb{R}^{[m,n]}$ be a strong $\M$-tensor. If $0<\omega_1<\omega_2\leq1$ and
$\alpha_{i},{\beta}_{j}\in[0,1],~i=1,...,n-s,~j=k+1,...,n$, then
\[\rho(\mathcal H_{\boldsymbol{\alpha}\boldsymbol{\beta}}(\omega_2))\leq\rho(\mathcal H_{\boldsymbol{\alpha}\boldsymbol{\beta}}(\omega_1))<1.\]
\end{prop}
\section{Numerical Examples}
In this section, we give some numerical examples to show the performance of the proposed algorithms. All tests were carried out in double precision with a \textsc{Matlab} code, when the computer specifications are Microsoft Windows 10 Intel(R), Core(TM)i7-7500U, CPU 2.70 GHz, with 8 GB of RAM. All used codes came from the \textsc{Matlab} tensor toolbox developed by Bader and Kolda \cite{A50,A40}.
We use PJ, PGS and PSOR to abbreviate the preconditioned Jacobi, Gauss-Seidel and SOR tensor splittings in \cite{A22}, \cite{A20} and \cite{A16}, respectively. In addition we use $P_{\boldsymbol{\alpha}\boldsymbol{\beta}}{\mathcal E_2\F_2},~P_{\boldsymbol{\alpha}\boldsymbol{\beta}}{\mathcal E_5\F_5},~ P_{\boldsymbol{\alpha}\boldsymbol{\beta}}{\M_2\N_2}, ~P_{\boldsymbol{\alpha}\boldsymbol{\beta}}{\M_5\N_5}$ and $P_{\boldsymbol{\alpha}\boldsymbol{\beta}}{SOR}$ to abbreviate the preconditioned Jacobi, Gauss-Seidel and SOR type splittings methods, respectively, that are proposed in this paper. We use Iter and Time for denote the number of iterations and CPU Times in seconds, respectively that need to reach the desired solution. The stoping criterion is $\norm{\textbf{r}_j}<10^{-12}$, where $\textbf{x}_0=\textbf{0}$, ${\textbf{r}_j=\textbf{b}-\A\textbf{x}_j^{m-1}}$ is the $j$-th iteration residual, the right hand side vector $\textbf{b}$ is $\textbf{1}=(1,...,1)^T$, if no other special illustration, and the maximum number of iterations is 2000. Also, we suppose that $\boldsymbol{\boldsymbol{\beta}}={\beta}\textbf{1}$ and $\boldsymbol{\alpha}=\alpha\textbf{1}$.\\
\textbf{Example 1.} Consider a strong $\M$-tensor $\A\in\Bbb{R}^{3\times3\times3}$, where
\[\tiny{\A(:,:,1)=\left(\begin{array}{*{20}c}
           1.00&-0.01&-0.02\\
           -0.02&-0.03&-0.04\\
           -0.04&-0.05&-0.06
          \end{array}\right),~
          \A(:,:,2)=\left(\begin{array}{*{20}c}
           -0.06&-0.07&-0.08\\
           -0.08&1.00&-0.09\\
           -0.01&-0.02&-0.03
          \end{array}\right),~
          \A(:,:,3)=\left(\begin{array}{*{20}c}
           -0.03&-0.04&-0.05\\
           -0.05&-0.06&-0.07\\
           -0.07&-0.08&1.00
          \end{array}\right).}\]
We compare the mentioned methods, where the parameter $\omega$ in the SOR method is chosen 1.2. We take ${\alpha}={\beta}$ in the interval $[0,10]$ with the step size 0.5 and $s,k=2$. The comparison results are shown in Table \ref{T1}.
\begin{table}[h]\small{
\caption{Iteration number (Iter) and CPU Time (Time) for Example 1.
\label{T1}}
\begin{center}
\begin{tabular}{|c|cccccc|}
\hline
&PJ&PGS&PSOR&$P_{\boldsymbol{\alpha}\boldsymbol{\beta}}{\mathcal E_2\F_2}$&$P_{\boldsymbol{\alpha}\boldsymbol{\beta}}{\M_2\N_2}$ &$P_{\boldsymbol{\alpha}\boldsymbol{\beta}}{SOR}$\\
\hline
${\alpha}$&Iter~~Time&Iter~~Time&Iter~~Time&Iter~~Time&Iter~~Time& Iter~~Time\\
\hline
0.0&51~~0.0066&50~~0.0065&39~~0.0100&51~~0.0048&50~~0.0055&39~~0.0050\\
0.5&51~~0.0046&49~~0.0041&39~~0.0044&50~~0.0018&49~~0.0030&38~~0.0018\\
1.0&50~~0.0044&47~~0.0032&39~~0.0033&49~~0.0028&48~~0.0019&37~~0.0014\\
1.5&50~~0.0039&46~~0.0040&39~~0.0024&48~~0.0017&47~~0.0018&36~~0.0017\\
2.0&50~~0.0040&45~~0.0020&39~~0.0021&46~~0.0017&46~~0.0017&35~~0.0012\\
2.5&49~~0.0040&44~~0.0020&39~~0.0021&45~~0.0017&45~~0.0017&35~~0.0016\\
3.0&49~~0.0030&43~~0.0024&39~~0.0020&44~~0.0016&44~~0.0016&34~~0.0012\\
3.5&48~~0.0028&43~~0.0019&39~~0.0019&43~~0.0015&43~~0.0013&33~~0.0011\\
4.0&48~~0.0032&42~~0.0023&39~~0.0027&42~~0.0022&42~~0.0013&32~~0.0012\\
4.5&47~~0.0021&43~~0.0021&39~~0.0019&41~~0.0015&41~~0.0012&31~~0.0013\\
5.0&47~~0.0023&43~~0.0023&39~~0.0028&40~~0.0016&40~~0.0012&30~~0.0010\\
5.5&46~~0.0023&44~~0.0024&39~~0.0021&39~~0.0021&39~~0.0014&29~~0.0023\\
6.0&46~~0.0026&44~~0.0022&39~~0.0019&38~~0.0022&38~~0.0014&28~~0.0016\\
6.5&45~~0.0024&45~~0.0022&39~~0.0020&37~~0.0016&36~~0.0019&27~~0.0012\\
7.0&44~~0.0022&46~~0.0024&39~~0.0020&35~~0.0020&35~~0.0014&26~~0.0022\\
7.5&44~~0.0027&48~~0.0024&39~~0.0019&34~~0.0017&33~~0.0013&24~~0.0012\\
8.0&43~~0.0023&49~~0.0025&39~~0.0021&\textbf{29~~0.0014}&31~~0.0013&22~~0.0009\\
8.5&43~~0.0029&50~~0.0027&39~~0.0021&31~~0.0016&28~~0.0012&\textbf{21~~0.0008}\\
9.0&42~~0.0026&52~~0.0028&39~~0.0022&33~~0.0018&\textbf{27~~0.0010}&23~~0.0013\\
9.5&42~~0.0024&53~~0.0028&39~~0.0022&34~~0.0017&31~~0.0014&25~~0.0014\\
10.0&41~~0.0024&55~~0.0028&39~~0.0022&34~~0.0017&32~~0.0014&28~~0.0017\\
\hline
\end{tabular}
\end{center}}
\end{table}
In addition, we take the $\omega$ in the interval $[0.5,1.8]$ with the step size 0.1 and obtain the solution by using the proposed preconditioned SOR method for $\alpha,{\beta}=0,5$ and $s,k=1,2$. We depicted the results in Table \ref{T2}, $P_{\boldsymbol{\alpha}}(s)=D+S_{\boldsymbol{\alpha}}(s)$ and $P_{\boldsymbol{\beta}}(k)=D+K_{\boldsymbol{\beta}}(k)$.
\begin{table}[h]\small{
\caption{Iteration numbers (Iter) and CPU Times (Time) for the preconditioned SOR type method.
\label{T2}}
\begin{center}
\begin{tabular}{|c|cccccc|}
\hline
&$P_{\boldsymbol{\alpha}=\textbf{5}}(1)$&$P_{\boldsymbol{\alpha}=\textbf{5}}(2)$ &$P_{\textbf{5,5}}(1,1)$ &$P_{\textbf{5,5}}(1,2)$&$P_{\textbf{5,5}}(2,1)$ &$P_{\textbf{5,5}}(2,2)$\\
\hline
$\omega$&Iter~~Time&Iter~~Time&Iter~~Time& Iter~~Time&Iter~~Time&Iter~~Time\\
\hline
0.5&103~~0.0206&95~~0.0106&105~~0.0119&\textbf{94~~0.0133}&97~~0.0128&96~~0.0148\\
0.6&83~~0.0089&{77~~0.0025}&85~~0.0036&\textbf{76~~0.0026}&79~~0.0027&77~~0.0027\\
0.7&69~~0.0018&{63~~0.0018}&70~~0.0020&63~~0.0018&\textbf{60~~0.0019}&64~~0.0024\\
0.8&58~~0.0017&\textbf{53~~0.0015}&60~~0.0022&53~~0.0019&56~~0.0017&54~~0.0016\\
0.9&50~~0.0014&\textbf{45~~0.0014}&51~~0.0018&46~~0.0014&48~~0.0014&46~~0.0016\\
1.0&43~~0.0014&\textbf{39~~0.0012}&44~~0.0013&40~~0.0012&42~~0.0014&40~~0.0012\\
1.1&37~~0.0012&\textbf{34~~0.0007}&39~~0.0007&35~~0.0010&37~~0.0011&35~~0.0010\\
1.2&33~~0.0009&\textbf{29~~0.0007}&34~~0.0010&30~~0.0009&32~~0.0019&30~~0.0008\\
1.3&29~~0.0008&\textbf{25~~0.0007}&30~~0.0008&26~~0.0007&29~~0.0008&26~~0.0008\\
1.4&31~~0.0009&\textbf{28~~0.0008}&30~~0.0009&29~~0.0008&33~~0.0011 &\textbf{28~~0.0010}\\
1.5&40~~0.0010&\textbf{35~~0.0009}&39~~0.0013&37~~0.0010&44~~0.0011 &\textbf{35~~0.0010}\\
1.6&53~~0.0014&\textbf{45~~0.0012}&51~~0.0013&47~~0.0012&59~~0.0017&46~~0.0019\\
1.7&71~~0.0028&\textbf{60~~0.0022}&70~~0.0020&64~~0.0017&81~~0.0013&61~~0.0010\\
1.8&105~~0.0015&\textbf{84~~0.0012}&101~~0.0014&92~~0.0013&135~~0.0020&85~~0.0019\\
\hline
\end{tabular}
\end{center}}
\end{table}
\newline
From Table \ref{T1}, we find that all the preconditioned methods perform better in CPU Times and iteration numbers than the ones with unpreconditioned ($\alpha=\beta=0$). Also, the proposed preconditioned schemes of Jacobi, Gauss-Seidel and SOR methods are all better than the corresponding ones that are considered in this paper when the parameters $\alpha$ and ${\beta}$ can be taken suitably. The best answers in terms of CPU times and iteration numbers have bolded in Table \ref{T1}. From Table \ref{T2} and for every choice of $\omega$, we see that in most cases when $\beta=0$ and $s=2$, the best answers in terms of CPU Times and iteration numbers are obtained which are showed in bolded numbers.
\newline
\textbf{Example 2.} Let $\B\in\Bbb{R}^{[3,n]}$ be a nonnegative tensor with $M(\B)$=hilb(n,n), where hilb is the function of \textsc{Matlab}, for $i=2,3,...,n,b_{ii-1i}=b_{iii-1}=b_{ii+1i}=b_{iii+1}=\tfrac{1}{3}$ and other entries are zeros. Let $\A=n^2\I-0.01\B$. We take $\alpha={\beta}=1,~s=k=n-1$ for $P_{\boldsymbol{\alpha}\boldsymbol{\beta}}\mathcal E_1\F_1,P_{\boldsymbol{\alpha}\boldsymbol{\beta}}\mathcal E_2\F_2,P_{\boldsymbol{\alpha}\boldsymbol{\beta}}SOR$. Also we obtained experimentally the optimal parameter $\omega$ in the interval $[0,2]$. The numerical results are reported in Table \ref{T3} which illustrate that the proposed preconditioned methods perform better in CPU times than the ones with the others.
\begin{table}[h]\small{
\caption{Iteration number (Iter) and CPU Time (Time) for Example 2.
\label{T3}}
\begin{center}
\begin{tabular}{|c|cccccccc|}
\hline
&PJ&PGS&PSOR&$P_{\boldsymbol{\alpha}\boldsymbol{\beta}}{\mathcal E_2\F_2}$&$P_{\boldsymbol{\alpha}\boldsymbol{\beta}}{\mathcal E_5\F_5}$& $P_{\boldsymbol{\alpha}\boldsymbol{\beta}}{\M_2\N_2}$& $P_{\boldsymbol{\alpha}\boldsymbol{\beta}}{\M_5\N_5}$ &$P_{\boldsymbol{\alpha}\boldsymbol{\beta}}{SOR}$\\
\hline
$n$&Iter~~Time&Iter~~Time&Iter~~Time&Iter~~Time&Iter~~Time& Iter~~Time&Iter~~Time&Iter~~Time\\
\hline
30&4~~0.0360&5~~0.0169&3~~0.0228&\textbf{3~~ 0.0151}&3~~0.0242 &3~~0.0196&
3~~0.0230&3~~0.0240\\
40&4~~0.0597&5~~0.0273&3~~0.0299&\textbf{3~~0.0174 }&3~~0.0327&3~~0.0233&
3~~0.0276&3~~0.0243\\
50&4~~0.0739&5~~0.0293&3~~0.0339&\textbf{3~~0.0201}&3~~0.0362&3~~0.0257&
3~~0.0306&3~~0.0288\\
60&4~~0.0608&5~~0.0443&3~~0.0483&\textbf{3~~0.0292}&3~~0.0499&3~~0.0417&
3~~0.0550&3~~0.0425\\
70&4~~0.0730&5~~0.0558&3~~0.0645 &\textbf{3~~0.0416}&3~~0.0787&3~~0.0575&
3~~0.0604&3~~0.0657\\
80&4~~0.0985&5~~0.0829&3~~0.0918&\textbf{3~~0.0497}&3~~0.1008&3~~0.0733&
3~~0.0829&3~~0.0929\\
90&4~~0.1173&5~~0.2064&3~~0.1300&\textbf{3~~0.0687}&3~~0.1329&3~~0.0859&
3~~0.0977&3~~0.1179\\
100&4~~0.1480&5~~0.1406&3~~0.1735&\textbf{3~~0.1173}&3~~0.2311&3~~0.1264&
3~~0.1373&3~~0.1449\\
110&4~~0.2133&5~~0.3544&3~~0.1968&\textbf{3~~0.1314}&3~~0.2432&3~~0.1525&
3~~0.1618&3~~0.1790\\
120&4~~0.2250&5~~0.2073&3~~0.2406&\textbf{3~~0.1464}&3~~0.2773&3~~0.1978&
3~~0.1967 &3~~0.2544\\
\hline
\end{tabular}
\end{center}}
\end{table}
\newline
From Table \ref{T3}, we find that when $n$ increases, the CPU Times for obtaining the appropriate answer increase. Also, if the parameters $\alpha,\beta$ and $\omega$ can be taken suitably, the proposed preconditioned schemes of Jacobi, Gauss-Seidel and SOR methods are all better than the corresponding ones that are considered in this paper. The best answers in terms of CPU Times and iteration numbers for every $n$ have bolded in Table \ref{T3}, where shows that the proposed second scheme of the preconditioned Jacobi method is the best.
\newline
\textbf{Example 3.} Let $\B\in\Bbb{R}^{[3,10]}$ be a nonnegative tensor and $b_{i_1i_2i_3}=|\tan(i_1+i_2+i_3)|$. It is not difficult (\cite{R23}) to see that $\rho(\B)\approx1450$, thus $\A=2000\I-\B$ is a strong $\M$-tensor. For mentioned methods, we obtained experimentally the optimal parameter $\omega$ in the interval $[1,2]$, the values of ${\alpha},{\beta}$ from 0 to 30 and $s,k=1$. The numerical results are reported in Table \ref{T4}. We use $\dagger$ to indicate that there was no convergence up to 2000 iterations. Table \ref{T4} illustrates that the proposed preconditioned methods perform better in CPU times than the ones with the others.
\begin{table}[h]{
\caption{Iteration number (Iter) and CPU Time (Time) for Example 3.
\label{T4}}
\begin{center}
\begin{tabular}{|cc|cccccc|}
\hline
&&PJ&PGS&PSOR&$P_{\boldsymbol{\alpha} \boldsymbol{\beta}} {\mathcal E_2\F_2}$&$P_{\boldsymbol{\alpha}\boldsymbol{\beta}}{\M_2\N_2}$ &$P_{\boldsymbol{\alpha}\boldsymbol{\beta}}{SOR}$\\
\hline
$\alpha$&$\beta$&Iter~~Time&Iter~~Time& Iter~~Time&Iter~~Time&Iter~~Time&Iter~~Time\\
\hline
0&0&91~~0.0308&87~~0.0181&69~~0.0302&91~~0.0168&
87~~0.0186&69~~0.0191\\
0.5&0.5&$\dagger$~~0.1253&$\dagger$~~0.0831&$\dagger$~~0.0820&90~~0.0143&
86~~0.0184&68~~0.0158\\
1&1&$\dagger$~~0.1253&$\dagger$~~0.0831&$\dagger$~~0.0820&89~~0.0147&
85~~0.0195&67~~0.0227\\
2&2&$\dagger$~~0.0960&$\dagger$~~0.0596&$\dagger$~~0.0853&87~~0.0149&
83~~0.0241&65~~0.0157\\
3&2&$\dagger$~~0.1025&$\dagger$~~0.0807&$\dagger$~~0.0949&85~~0.0169&
81~~0.0167&64~~0.0152\\
4&2&$\dagger$~~0.1051&$\dagger$~~0.0740&$\dagger$~~0.0788&83~~0.0145&
79~~0.0204&62~~0.0148\\
5&5&$\dagger$~~0.1023&$\dagger$~~0.0745&$\dagger$~~0.0849&81~~0.0167&
77~~0.0201&60~~0.0160\\
7&5&$\dagger$~~0.0953&$\dagger$~~0.0999&$\dagger$~~0.0831&77~~0.0144&
73~~0.0185&57~~0.0160\\
9&5&$\dagger$~~0.0835&$\dagger$~~0.0822&$\dagger$~~0.0989&73~~0.0142&
69~~0.0164&54~~0.0150\\
10&8&$\dagger$~~0.1347&$\dagger$~~0.1228&$\dagger$~~0.1262&71~~0.0155&
67~~0.0193&52~~0.0146\\
12&10&$\dagger$~~0.0975&$\dagger$~~0.0848&$\dagger$~~0.0760&67~~0.0159&
64~~0.0159&49~~0.0145\\
15&12&$\dagger$~~0.0934&$\dagger$~~0.0875&$\dagger$~~0.0987&61~~0.0134&
58~~0.0167&44~~0.0153\\
18&10&$\dagger$~~0.0912&$\dagger$~~0.0985&$\dagger$~~0.1044&55~~0.0166&
52~~0.0186&39~~0.0150\\
20&15&$\dagger$~~0.0987&$\dagger$~~0.0924&$\dagger$~~0.0901&49~~0.0138&
48~~0.0157&38~~0.0137\\
20&20&$\dagger$~~0.0912&$\dagger$~~0.0914&$\dagger$~~0.0926&49~~0.0144&
48~~0.0159&38~~0.0182\\
25&20&$\dagger$~~0.0989&$\dagger$~~0.0932&$\dagger$~~0.0911&\textbf{40~~0.0138}&
\textbf{40~~0.0151}&38~~0.0139\\
25&25&$\dagger$~~0.0999&$\dagger$~~0.0924&$\dagger$~~0.0937&42~~0.0139&
\textbf{40~~0.0153}&\textbf{37~~0.0137}\\
30&20&$\dagger$~~0.0974&$\dagger$~~0.0978&$\dagger$~~0.0945&47~~0.0140&
48~~0.0166&45~~0.0140\\
30&25&$\dagger$~~0.0910&$\dagger$~~0.0934&$\dagger$~~0.0922&48~~0.0139&
48~~0.0170&46~~0.0140\\
30&30&$\dagger$~~0.0900&$\dagger$~~0.0944&$\dagger$~~0.0891&48~~0.0141&
49~~0.0158&46~~0.0140\\
\hline
\end{tabular}
\end{center}}
\end{table}
\newline
As we see in Table \ref{T4}, without preconditioners ($\alpha=\beta=0$), the proposed preconditioned schemes of Jacobi, Gauss-Seidel and SOR methods obtained the same answers with the corresponding ones that are considered in this paper. When the parameters $\alpha$ and $\beta$ are considered as nonzero, we see that the PJ, PGS and PSOR methods are not convergent, but the proposed methods are convergent and improve the iteration numbers and CPU Times concerning unpreconditioned. The best answers in the iteration numbers and CPU Times are bolded in Table \ref{T4}.
\section{Concluding remarks}
In this paper, we proposed new types of flexible and fast preconditioners tensor splitting methods for solving multilinear system $\A\textbf{x}^{m-1}=\textbf{b}$, when $\A$ is a strong $\M$-tensor. Some properties of convergent theorems about preconditioned Jacobi, Gauss-Seidel and SOR type iterative methods are obtained. Numerical examples are given to show the efficiency and superiority of the proposed methods.

\end{document}